\documentclass[]{article}

\usepackage{amsmath,amsthm,amsfonts}
\usepackage{geometry}
\usepackage{hyperref}

\newcommand{\boldone}{\mathbf{1}}
\newcommand{\ee}{\mathrm{e}}
\newcommand{\BE}{\mathbb{E}}
\newcommand{\BP}{\mathbb{P}}

\newcommand{\floor}[1]{\lfloor#1\rfloor}
\newcommand{\deq}{\mathrel{\overset{d}{=}}}

\newtheorem*{theorem*}{Theorem}

\newcommand{\MR}[1]{\href{http://www.ams.org/mathscinet-getitem?mr=#1}{MR#1}}

\title{A simple proof\\ of the L\'evy--Khintchine formula for subordinators}
\author{Yuri Yakubovich\thanks{St. Petersburg State University, 7/9 Universitetskaya nab., St. Petersburg,
		199034 Russia.\newline E-mail: \href{mailto:y.yakubovich@spbu.ru}{y.yakubovich@spbu.ru}\newline
		Partially supported by the Russian Foundation
		for Basic Research, project no.~20-01-00646~A}}

\begin{document}

\maketitle

\begin{abstract}
We present a relatively simple and mostly elementary proof of the L\'evy--Khintchine formula for subordinators. 
The main idea is to study the Poisson process time-changed by the subordinator. The technical tools used are conditional
expectations, probability generating function and convergence of discrete random variables.
\end{abstract}

\section{Introduction}
A subordinator is a real-valued stochastic process $S(t)$, $t\ge 0$ 
which starts at 0 (that is $S(0)=0$), does not decrease ($S(s)\le S(t)$ for any $0\le s\le t$), is right-continuous ($\lim_{t\downarrow s}S(t)=S(s)$ 
at every point $s$),
and has stationary increments (the distribution of $S(t)-S(s)$ depends only on the difference $t-s$, for any $0\le s\le t$)
 which are independent (for every $n$ and any $0=t_0\le t_1\le \dots\le t_n$ the increments
 $S(t_1)-S(t_0), S(t_2)-S(t_1),\dots,S(t_n)-S(t_{n-1})$ are independent random variables).\footnote{%
To make our presentation as simple as possible we avoid using notions of the theory of stochastic processes,
in particular filtrations which allow to formulate the assumptions shorter, and confine ourselves to notions of the elementary probability theory.} 
By a stochastic process we mean a collection of random variables $S(t,\omega)$ (where $\omega\in\Omega$ for some
abstract probability space $(\Omega,\mathcal{F},\BP)$)  indexed by a real parameter $t\ge 0$. 
For a fixed $\omega$ it is a function of $t$, known as sample trajectory of a process, and this function has the stated 
properties almost surely (a.s.),\ that is for every $\omega\in\Omega'$ for 
some $\Omega'\subset\Omega$ such that $\BP(\Omega')=1$. We often refer to the parameter $t$ as time. The distribution of a process is usually
given in terms of joint distributions of its values in finitely many arbitrary times.  But for a process with stationary and independent
increments it obviously suffices to specify just one-dimensional distributions of $S(t)$, for any $t\ge 0$. Slightly less obvious 
is that actually it suffices to do it for some $t>0$.

The simplest (and not interesting) example of a subordinator is a non-random line $S(t)=\beta t$ for some $\beta\ge 0$
(here the distribution of $S(t)$ is the delta-measure at $\beta t$). We shall exclude this degenerate case and suppose that for
some (and hence for all) $t>0$ the random variable $S(t)$ has a non-degenerate distribution.   

Another simple but more interesting example is the standard Poisson process. It can be defined as follows.  Let 
$Y_1,Y_2,\dots$ be a sequence of independent identically distributed (i.i.d.)\ exponential random variables with mean 1.
Denote $T_n=Y_1+\dots+Y_n$ the sum of the first $n$ members of this sequence and let 
\begin{equation}\label{Poisson}
\Pi(t)=\sum_{n=1}^{\infty}\boldone_{\{T_n\le t\}}.
\end{equation}
Here and below $\boldone_A$ is the indicator of an event $A$, that is a random variable which takes value 1 on $A$ and 0 on its complement.
By the strong law of large numbers $T_n\to\infty $ as $n\to\infty$ a.s.,\ so the sum in \eqref{Poisson} is actually finite
with probability 1, and $\Pi(t)$ is well defined. By definition, $\Pi(0)=0$ and $\Pi$ does not decrease. The stationarity and independence
of increments follow easily from the memoryless property of the exponential distribution, that is $\BP(Y_1\ge t+s|Y_1\ge t)=\BP(Y_1\ge s)$
for any $s,t>0$.  The name Poisson process is connected with the fact that the distribution
of $\Pi(t)$ at any fixed time $t$ is Poisson: $\BP[\Pi(t)=k]=t^k \ee^{-t}/k!$ for all $k=0,1,\dots$.  
A generic Poisson process of \textit{rate} $\lambda>0$ is a linear time change of the standard one, that is the process $\Pi(\lambda t)$,
so standard means ``rate~1''.
We refer to \cite[Ch.~2.6]{Durrett} for an elementary exposition of these facts. 

The trajectories of any Poisson process are piecewise constant and have jumps at random moments of time $T_1,T_2,\dots$.  Each jump
is just a unit step. If we relax this requirement and allow steps to be random but i.i.d.\ we get an instance of a compound Poisson process.
In general, steps can be negative, but in order to obtain a subordinator we allow just positive real steps. Thus a compound Poisson 
subordinator of rate $\lambda$ is a process
\begin{equation}\label{compound}
C\!P(t)=\sum_{k=1}^{\Pi(\lambda t)}J_k,
\end{equation}
where $J_1,J_2,\dots$ are i.i.d.\ positive random variables.

The class of compound Poisson subordinators with some non-random linear drift $C\!P(t)+\beta t$, $\beta\ge 0$, is known as subordinators with finite
activity.  There is however a wide class of subordinators which are not compound Poisson processes, and exhibit quite different properties: 
they have infinitely many jumps in any finite interval, so their
trajectories are not piecewise constant or linear, but have jumps in a neighborhood of any point.  The complete characterization of possible laws of subordinators is well known
in terms of the Laplace transform of $S(t)$:
\begin{equation}\label{Levy=Khintchine}
\BE\bigl[\ee^{-uS(t)}\bigr]=\ee^{-t\Psi(u)},\qquad \Psi(u)=\beta u +\int_0^\infty (1-\ee^{-ux})\Lambda(dx),\qquad u\ge 0,
\end{equation} 
where $\beta\ge 0$ is called \textit{drift} coefficient, and $\Lambda$ is a measure on $(0,\infty)$ such that the integral converges for
some (and hence for all) $u>0$, known as the \textit{L\'evy measure}.  (Another common way to impose the same condition is to say that 
$\int_0^\infty\min\{x,1\}\Lambda(dx)<\infty$).  The function $\Psi(u)$ is called the \textit{Laplace exponent} of the subordinator. 
To be more precise, let us formulate the result.

\begin{theorem*}[L\'evy--Khintchine representation for subordinators]
Any function $\Psi$ given by~\eqref{Levy=Khintchine} is the Laplace exponent of some subordinator $S(t)$, $t\ge 0$.
Conversely, any subordinator $S(t)$, $t\ge 0$, admits a representation \eqref{Levy=Khintchine}, and the drift $\beta$ and the
L\'evy measure $\Lambda$ are uniquely determined.
\end{theorem*}

The first part of this claim is almost obvious.  Indeed, if \eqref{Levy=Khintchine} holds $S(0)$ is 0 a.s.\ and independence, 
positivity and stationarity of increments are compatible with \eqref{Levy=Khintchine}.  
Traditional ways to proof the converse go back to works of de Finetti, L\'evy,  Khintchine and Kolmogorov (we
refer to \cite{Mainardi} for a review of their pioneering work) and have analytical nature.
For a modern exposition see \cite{Bertoin,Sato}. The aim of this note is to provide an apparently new, 
more probabilistic proof for the L\'evy--Khintchine theorem
and to keep it as elementary as possible. This is achieved via random time change in the standard Poisson process $\Pi(s)$ by
setting $s=S(t)$, for independent $\Pi$ and $S$. This procedure introduced in \cite{Bochner} is known 
as a subordination (this explains the name) of the process $\Pi$ 
by the subordinator $S$; however, usually the subordinated process has more
complex structure than the Poisson process, and sometimes is dependent on the subordinator.

In the next section we present our proof.  Some remarks complete this note.

\section{A proof}
Let a subordinator $S(t)$, $t\ge 0$, be defined on a probability space $(\Omega,\mathcal{F},\BP)$ rich enough to allow an independent
standard Poisson process $\Pi(s)$, $s\ge 0$, to be defined on the same space.  Consider a composition $\Pi(S(t))$.  It is an integer-valued
non-decreasing random process.  It is easy to see that its increments are independent and stationary, hence it is also a subordinator.
Our first aim is to show that it is a compound Poisson process.  Let $t>0$.  Since $\Pi(S(t))$ is a.s.\ finite,
and by construction $\Pi(S(\tau))$ does not decrease for $0\le \tau\le t$ and admits just non-negative integer values, 
it has a.s.\ finite number $K(t)$ of jumps 
($K(t)=0$ if $\Pi(S(t))=0$), say at (random) times $0<\tau_1<\dots<\tau_{K(t)}<t$.  Let $\zeta=\min_{i=1,\dots,K(t)+1}\{\tau_i-\tau_{i-1}\}$ 
(where $\tau_0=0$ and $\tau_{K(t)+1}=t$) be the length of the shortest interval between these jump times.

Introduce a notation $\floor{s}$ for the floor function, that is $\floor{s}$ is the maximal integer not greater than $s$. 
For a natural number $n$ divide the interval $[0,t]$ on $\floor{nt}$ intervals 
$\Delta_{n,i}:=\bigl[\tfrac{i-1}{n},\tfrac{i}{n}\bigr]$ of length $\tfrac1n$ and the last,
shorter interval $\Delta_n^*:=\bigl[\frac{\floor{nt}}{n},t\bigr]$.  
For each $n$ introduce two sequences of random variables,
indicators $I_{n,i}=\boldone_{\{\Pi(S(i/n)) > \Pi(S((i-1)/n))\}}$ of the events that $\Pi(S(\cdot))$ jumps on the $i$th interval ($i=1,\dots,\floor{nt}$),
and jumps $J_{n,k}=\Pi(S(\tfrac{\floor{n\tau_k}+1}{n}))-\Pi(S(\tfrac{\floor{n\tau_k}}{n}))$ of $\Pi(S(\cdot))$ 
on the interval containing $\tau_k$ ($k=1,\dots,K(t)$).
For any fixed $n$ the indicators $I_{n,i}$ are i.i.d.\ by stationarity and independence of increments of $\Pi(S(\cdot))$. The jumps
$J_{n,k}$ also are identically distributed but in general not independent: if, say, $\tau_k$ and $\tau_{k+1}$
fall in the same interval $\Delta_{n,i}$ then $J_{n,k}=J_{n,{k+1}}$. As $n$ grows, however, they eventually (at least for $n\ge 1/\zeta$)
start measuring jumps on different intervals and hence independent. So their a.s.\ limits $J_k:=\lim_{n\to\infty}J_{n,k}=\Pi(S(\tau_k))-\Pi(S(\tau_k-0))$
are i.i.d.\footnote{This argument actually contains some cheating: since $\zeta=\zeta(\omega)$ can be arbitrary close to 0 one can not 
claim that the intervals containing, say, $\tau_k$ and $\tau_{k+1}$ become different 
for sufficiently large $n$ for almost all $\omega\in\Omega$,  so $J_{n,k}$ and $J_{n,k+1}$ are not independent for any $n$. It also does not 
look obvious that $J_{n,k}$ and $J_{n,k+1}$ are conditionally independent given $\zeta<n$ because the last event affects the whole
trajectory of $\Pi(S(\cdot))$.  However since $\BP[\zeta<n]\to 1$ as $n\to\infty$, 
the conditional probability differs very little from the unconditional one for large
$n$, and $J_{n,k}$ and $J_{n,k+1}$ are asymptotically independent and independent in the limit.}

Let $K_n(t)=\sum_{i=1}^{\floor{nt}}$ be
the number of intervals $\Delta_{n,i}$ (ignoring  $\Delta_n^*$) on which $\Pi(S(\cdot))$ has at least one jump.  For $n>1/\zeta$ the subordinator
$\Pi(S(\cdot))$ jumps not more than once on each of these intervals, and does not jump on $\Delta_n^*$.
Hence $K_n(t,\omega)=K(t,\omega)$ for $\omega\in\Omega$ such that $n>1/\zeta(\omega)$, 
and $K_n(t)\to K(t)$ a.s.\ as $n\to \infty$.  Consequently, the distribution of $K_n(t)$ converges to the
distribution of $K(t)$.

Let us find the distribution of $K_n(t)$ directly. Denote $q_n=\BP[I_{n,i}=0]$.  The intersection of the first $n$ such events is the event $\cap_{i=1}^n\{I_{n,i}=0\}=\bigl\{\Pi(S(1))=0\bigr\}$
which has probability
\begin{equation}\label{qnn}
q_n^n=\BP\bigl[\Pi(S(1))=0\bigr]=\BE\bigl[ \boldone_{\{\Pi(S(1))=0\}}\bigr]
=\BE\bigl[\BE \bigl(\boldone_{\{\Pi(S(1))=0\}}\big|S(1)\bigr)\bigr]=\BE\bigl[\ee^{-S(1)}\bigr].
\end{equation}
Denote this probability by $\ee^{-\psi}$, for some $\psi\ge 0$.  Then $q_n=\ee^{-\psi/n}$ and $K_n(t)$ has a binomial distribution 
with the probability generating function
\[
\BE\bigl[z^{K_n(t)}\bigr]=(q_n+z(1-q_n))^{\floor{nt}}=(1-(1-z)(1-q_n))^{\floor{nt}}.
\]
Since $1-q_n=\psi/n+O(n^{-2})$ this converges, as $n\to\infty$, to
\[
\lim_{n\to\infty}\BE\bigl[z^{K_n(t)}\bigr]=\ee^{-t\psi (1-z)}
\]
which is the probability generating function of the Poisson distribution with mean $t\psi$ implying (see \cite[Ch.~XI.6]{Feller}) 
that $K(t)$ has the same distribution.  
Combining the a.s.\ convergences together gives the representation of $\Pi(S(t))$ as a sum of its jumps at discontinuities
\begin{equation}\label{Pi(S)-compound}
\Pi(S(t))=\sum_{k=1}^{K(t)} J_k
\end{equation}
so it is indeed a compound Poisson process.

Consider now jumps $J_{n,k}$. They all have the distribution of the $\Pi(S(\tfrac1n))$ conditioned to have at least one jump on $[0,1/n]$:
$J_{n,k}\deq\bigl(\Pi(S(\tfrac1n))\bigm|\Pi(S(\tfrac1n))>0\bigr)$.  For a fixed $s$ the probability generating function of such conditional
distribution is
\[
\BE\bigl[z^{\Pi(s)}\big|\Pi(s)>0\bigr]=\sum_{k=1}^\infty z^k\frac{s^k \ee^{-s}}{k!(1-\ee^{-s})}=\frac{\ee^{(z-1)s}-\ee^{-s}}{1-\ee^{-s}}
\]
where the expectation is taken with respect to the conditional probability given $\Pi(s)>0$. As a function of $s$ it does not increase for any fixed
$z\in[0,1]$, so
\[
\BE\bigl[z^{\Pi(S(1/n))}\big|\Pi(S(1/n))>0\bigr]
=\BE\Bigl[\frac{\ee^{(z-1)S(1/n)}-\ee^{-S(1/n)}}{1-\ee^{-S(1/n)}}\Bigr]
\]
increases in $n$ for a fixed $z\in[0,1]$ (because $S(\tfrac{1}{n+1})\le S(\tfrac1n)$ a.s.)\ and has a limit
\[
\lim_{n\to\infty} \BE\bigl[z^{\Pi(S(1/n))}\big|\Pi(S(1/n))>0\bigr] = :f(z).
\]
This limit is the probability generating function of any jump $J_k=\Pi(S(\tau_k))-\Pi(S(\tau_k-0))\deq\Pi(S(\tau_k)-S(\tau_k-0))$ and can
be represented as a similar expectation. We write it as the Stieltjes integral with respect to the distribution function
$F(x)=\BP[S(\tau_k)-S(\tau_k-0)\le x]$.  Note that $F(0)$ is equal to the probability
that $\Pi(\cdot)$ jumps at a point $s$ such that $S(\tau_k)=S(\tau_k-0)=s$ which can be positive.  Hence
\begin{equation}\label{f(x)}
f(z)=\BE[z^{J_k}]=\int_{[0,\infty)}\frac{\ee^{(z-1)x}-\ee^{-x}}{1-\ee^{-x}}dF(x).
\end{equation}

According to \eqref{Pi(S)-compound} we can write the probability generating function of $\Pi(S(t))$ as
\begin{align}\label{genfunc-int}
\notag
\BE\bigl[z^{\Pi(S(t))}\bigr]& = \sum_{k=0}^\infty \frac{(\psi t f(z))^k}{k!}\ee^{-\psi t}
= \ee^{\psi t (f(z)-1)}\\
\notag
&=\exp\Bigl(\psi t \int_{[0,\infty)}\Bigl(\frac{\ee^{(z-1)x}-\ee^{-x}}{1-\ee^{-x}}-1\Bigr)dF(x)
\Bigr)\\
&=\exp\Bigl(\psi t F(\{0\})(z-1)+\psi t \int_{(0,\infty)}\frac{\ee^{(z-1)x}-1}{1-\ee^{-x}}dF(x)\Bigr)
\end{align}
where in the second line we used $1=\int_{[0,\infty)}dF(x)$ because $dF(x)$ is the probability distribution, 
and in last line we have separated the atom $F(\{0\})$ at 0 of the measure $dF(x)$.  On the other hand, 
\begin{equation}\label{genfunc-exp}
\BE\bigl[z^{\Pi(S(t))}\bigr]=\BE\bigl[\BE\bigr(z^{\Pi(S(t))}\big|S(t)\bigr)\bigr]=\BE\bigl[\ee^{(z-1)S(t)}\bigr].  
\end{equation}
Let $\Lambda(dx)$ be the measure on $(0,\infty)$ with the density 
$\psi/(1-\ee^{-x})$ relative to $dF(x)$ and $\beta=\psi F(\{0\})$.  Denoting $u=1-z$ one obtains that \eqref{Levy=Khintchine} 
holds for $u\in[0,1]$ by equating \eqref{genfunc-int} and \eqref{genfunc-exp}, and also for complex $u$ such that $|u-1|\le 1$.  
But the Laplace transform $\BE[\ee^{-u S(t)}]$ is the complex analytic function of $u$ at least in the  
domain $\operatorname{Re}u\ge 0$ and hence is uniquely defined by \eqref{Levy=Khintchine} in all this domain (see, e.g.,~\cite[Ch. II.6]{Widder}).

\section{Final remarks}

\paragraph{Killed subordinator}
Some authors (e.g.~\cite{Bertoin}) define subordinators as processes with values in $[0,\infty]$, with the special
value $+\infty$ which is never leaved once reached, so if $S(s)=+\infty$ then $S(t)=+\infty$ for all $t\ge s$.  It is 
well understood that if $\BP[S(s)=+\infty]>0$ for some $s<\infty$ then 
the process $S(\cdot)$ can be obtained from some a.s.\ finite subordinator $\tilde S(\cdot)$ 
by introducing an independent exponential random variable $Y$ with mean $1/\alpha$, for
some $\alpha>0$ called \textit{killing rate}, and ``killing'' the subordinator (that is, sending it to $+\infty$) at random time $Y$:
\[
S(t)=\begin{cases}
	\tilde S(t), & t<Y,\\
	+\infty, & t\ge Y.
\end{cases}
\]
With the convention
that $\ee^{-\infty}=0$ this leads to a modification of the Laplace exponent:
\[
\Psi(u)=\alpha + \beta u + \int_0^\infty (1-\ee^{-ux})\Lambda(dx),\qquad u\ge 0.
\]
In this note we considered a.s.\ finite subordinators (sometimes called \textit{strict}) for which $\alpha=0$.

\paragraph{The leading process}
As we have seen, the process $K(t)$ counting jumps of $\Pi(S(\cdot))$ but ignoring their amplitude is the Poisson process of rate $\psi$.
Comparing \eqref{qnn} and \eqref{Levy=Khintchine} yields $\psi=\Psi(1)$.  So the only characteristic of $S(\cdot)$ that affects how often $\Pi(S(\cdot))$ 
jumps is $\Psi(1)$.  All other characteristics of $S(\cdot)$ are reflected just in the jump distribution of $\Pi(S(\cdot))$.

\paragraph{The jumps of $\Pi(S(\cdot))$}
It is easy to find the L\'evy measure $M$ of the compound Poisson process $\Pi(S(\cdot))$: it has atoms $m_j$ at $j=1,2,\dots$ given by
\begin{equation}\label{mj}
m_1=\beta+\int_0^\infty x\ee^{-x}\Lambda(dx)\text{ and }m_j=\int_0^\infty \frac{x^j}{j!}\ee^{-x}\Lambda(dx)\text{ for }j\ge 2.
\end{equation}
After dividing by the total mass $M(0,\infty)=\psi=\Psi(1)$ the L\'evy measure $M$ becomes the probability distribution of jumps $J_k$,
and \eqref{mj} follows from \eqref{f(x)} by a change of notation; cf.\ eq.~(1.1) in~\cite{Orsingher}.
In terms of the Laplace exponent $\Psi(u)$ the probability generating function of $J_k$ is
\begin{align*}
f(z)=\BE[z^{J-k}]&=\frac{1}{\psi}\sum_{j=1} z^j m_j = \frac{z\beta}{\psi} +\frac1\psi\int_0^\infty \sum_{j=1}^\infty \frac{(zx)^j}{j!}\ee^{-x}\Lambda(dx)\\
&= \frac{(z\beta}{\psi} +\frac1\psi  \int_0^\infty \bigl(\ee^{(z-1)x}-\ee^{-x}\bigr)\Lambda(dx)\\
&= \frac\beta\psi-\frac{(1-z)\beta}{\psi} +\frac1\psi  \int_0^\infty \bigl(\ee^{-(1-z)x}-1-(\ee^{-x}-1)\bigr)\Lambda(dx)\\
&=\frac{\Psi(1)-\Psi(1-z)}{\Psi(1)}.
\end{align*}

\paragraph{The jumps of $S(\cdot)$ at $\tau_k$}  
The times $\tau_k$ can be interpreted as times when the subordinator $S(t)$ jumps over the points of the Poisson point process
of intensity 1.  The distribution of the jumps $S(\tau_k)-S(\tau_k=0)$ is $\tfrac{1}{\Psi(1)}\bigl((1-\ee^{-x})\Lambda(dx)+\beta\delta_0\bigr)$.
When the process $S(\cdot)$ is itself compound Poisson \eqref{compound} it is just the distribution of its jump $J_1$ conditioned to be greater than 
an independent exponential random variable with mean 1. For infinite activity subordinators one can not speak about individual jumps 
but in certain sense it is also true.

\paragraph{Related research}
Poisson processes with a random time change by a subordinator recently became objects of extensive studies since for some specific 
subordinators, in particular for stable subordinators, they appear in application as one of the possible ways to define 
a fractional analog of the Poisson process, see, for instance, \cite{Garra,Orsingher,Buchak} and references therein.  They present more deep 
analysis of these processes than ours, in particular they study distributions of the first passage times and of hitting probabilities.


\begin{thebibliography}{10}
\bibitem{Bertoin} J. Bertoin. \textit{Subordinators: examples and applications.} Lectures on probability theory and statistics (Saint-Flour, 1997), 1--91, Lecture Notes in Math., 1717, Springer, Berlin, 1999. \MR{1746300}

\bibitem{Bochner}
S. Bochner. Harmonic analysis and the theory of probability. University
of California Press, Berkeley, 1955.
\MR{0072370}

\bibitem{Buchak}
K. Buchak, L. Sakhno. Properties of Poisson processes directed by compound Poisson-Gamma subordinators, \textit{Modern Stoch. Theory Appl.} 
\textbf{5} (2018), no.~2, 167--189. \MR{3813090}

\bibitem{Durrett} 
R. Durrett. \textit{Probability: theory and examples.} Second edition. Duxbury Press, Belmont, CA, 1996. 
\MR{1609153}

\bibitem{Feller}
W. Feller. An introduction to probability theory and its applications. Vol. I. Third edition. John Wiley \& Sons, Inc., New York-London-Sydney, 1968. 
\MR{0228020}

\bibitem{Garra}
R. Garra, E. Orsingher, M. Scavino. Some probabilistic properties of fractional point processes. \textit{Stoch. Anal. Appl.} 
\textbf{35} (2017), no.~4, 701--718. \MR{3651139}


\bibitem{Mainardi}
F. Mainardi, S. Rogosin. The origin of infinitely divisible distributions: from de Finetti's problem to Lévy-Khintchine formula. \textit{Math. Methods Econ. Finance} \textbf{1} (2006), no. 1, 37--55. \MR{2317505}

\bibitem{Orsingher}
E. Orsingher, B. Toaldo. Counting processes with Bern\v{s}tein intertimes and random jumps. \textit{J. Appl. Probab.} 
\textbf{52} (2015), no.~4, 1028--1044. \MR{3439170}

\bibitem{Sato} K.-I. Sato. \textit{L\'evy processes and infinitely
divisible distributions.} Cambridge Univ. Press, 1999.
\MR{1739520}

\bibitem{Widder} D.~V. Widder. \textit{The Laplace transform.} Princeton Univ. Press, 1941.
\MR{0005923}

\end{thebibliography}
\end{document}